\documentclass[11pt,a4wide,oneside,reqno]{amsart}
\usepackage{amssymb,amsmath,url,a4wide}
\usepackage[utf8]{inputenc}
\usepackage[left=3cm,right=3cm,top=3cm,bottom=2cm]{geometry}
\usepackage{enumitem}

\newcommand{\citep}[1]{\cite{#1}}
\newcommand{\Sym}{\mathrm{Sym}}
\newcommand{\cW}{\mathfrak{W}}

\newcommand{\F}{\mathbb{F}}
\newcommand{\N}{\mathbb{N}}

\newcommand{\e}{\varepsilon}

\newcommand{\rv}{\vert}
\newcommand{\lv}{\vert}
\renewcommand{\mod}{\;\mathrm{mod}\;}

\setlist[enumerate]{label=(\roman*)}

\usepackage[breaklinks=true]{hyperref}
\usepackage{cleveref}

\numberwithin{equation}{section}
\newtheorem{theorem}{Theorem}[section]
\newtheorem{claim}[theorem]{Claim}
\newtheorem{lemma}[theorem]{Lemma}

\newtheorem{proposition}[theorem]{Proposition}

\theoremstyle{definition}

\newtheorem{example}[theorem]{Example}

\usepackage{thmtools}
\declaretheoremstyle[%
  spaceabove=0pt,%
  spacebelow=0pt,%
  headfont=\normalfont\itshape,%
  postheadspace=1em,%
  qed=\qedsymbol%
]{mystyle}

\usepackage[url=false,doi=true,maxnames=4,giveninits=true,isbn=false]{biblatex}
\usepackage{xpatch}
\xpatchbibmacro{volume+number+eid}{%
  \setunit*{\adddot}%
}{%
}{}{}
\DeclareFieldFormat[article]{number}{\mkbibparens{#1}}
\renewbibmacro{in:}{%
  \ifentrytype{article}{}{\printtext{\bibstring{in}\intitlepunct}}}
\bibliography{biblio}

\title{The structure of large sum-free sets in $\mathbb{F}_p^n$}

\author{Leo Versteegen}
\thanks{The author is grateful to be funded by Trinity College of the University of Cambridge through the Trinity External Researcher Studentship.}
\address{Department of Pure Mathematics and Mathematical Statistics, Centre for Mathematical Sciences, Wilberforce Road, Cambridge CB3 0WB, United Kingdom}
\email{lvv23@dpmms.cam.ac.uk}

\begin{document}
\maketitle
\begin{abstract}
A set $A\subset \mathbb{F}_p^n$ is sum-free if $A+A$ does not intersect $A$. If $p\equiv 2 \mod 3$, the maximal size of a sum-free in $\mathbb{F}_p^n$ is known to be $(p^n+p^{n-1})/3$. We show that if a sum-free set $A\subset \mathbb{F}_p^n$ has size at least $p^n/3-p^{n-1}/6+p^{n-2}$, then there exists subspace $V<\mathbb{F}_p^n$ of co-dimension 1 such that $A$ is contained in $(p+1)/3$ cosets of $V$. For $p=5$ specifically, we show the stronger result that every sum-free set of size larger than $1.2\cdot 5^{n-1}$ has this property, thus improving on a recent theorem of Lev.
\end{abstract}
\section{Introduction}
Let $G$ be a finite Abelian group. Given subsets $A,B\subset G$, we write $A+B$ for their \emph{sumset} $\{a+b:a\in A, b\in B\}$. A set $A\subset G$ is said to be \emph{sum-free} if $(A+A)\cap A = \emptyset$, i.e., if there are no $a,b,c\in A$ such that $a+b=c$. As far back as 1969, Yap \cite{yap69} investigated the size of a largest sum-free set in $G$, denoted by $\lambda(G)$. Although there was some early progress on this problem for large classes of groups \cite{diananda_yap69, rhetmulla_street70, yap72,yap75}, the question was only fully resolved by Green and Rusza \cite{green_rusza05}, who determined $\lambda(G)$ for all finite Abelian groups. We confine ourselves to quoting the end-result of over 35 years of research and refer the reader to \cite{green_rusza05} for a more detailed account.

\begin{theorem}\label{thm:max-sizes}
Writing $N=\lv G\rv$, $\lambda(G)$ is determined as follows.
\begin{itemize}
\item If $N$ is divisible by a prime $p\equiv 2\mod 3$, then $\lambda(G)=N(p+1)/3p$, where $p$ is the smallest such prime.
\item If $N$ is not divisible by a prime $p\equiv 2\mod 3$ but $3\vert n$, then $\lambda(G)=N/3$.
\item If $N$ is divisible only by primes $p\equiv 1\mod 3$, then $\lambda(G)=(1-1/m)N/3$, where $m$ is the largest order of any element in $G$.
\end{itemize}
\end{theorem}

Long before \cite{green_rusza05}, there had already been interest in a related line of inquiry that asks what structure a large sum-free set may have. Already in \cite{yap69}, Yap observes that if the smallest prime divisor $p$ of $\lv G\rv$ is $2\mod 3$ and $A$ is sum-free with $\lv A\rv = \lambda(G) = \lv G\rv(p+1)/3p$, then $A$ is a union of cosets of a subgroup $H\subset G$ and $A/H$ is in arithemtic progression. As another example, motivated by applications in coding theory, Davydov and Tombak \cite{davydov_tombak89} and independently Clark and Pedersen \cite{clark_pedersen92} proved the following result.

\begin{theorem}\label{thm:p=2}
If $A\subset \F_2^n$ is sum-free and $\lv A\rv > 5\cdot 2^{n-4}$, then there exists a proper subspace $V<\F_2^n$ such that $A$ is contained in a coset of $V$.
\end{theorem}

The set $A=\{(1,0,0,0),(0,1,0,0),(0,0,1,0),(0,0,0,1),(1,1,1,1)\}$ shows that the bound $5\cdot 2^{n-4}$ is sharp. Later, Lev \cite{lev05} proved an analogous statement for $p=3$. 

In this paper, we are concerned with the structure of large sum-free subsets of $\F_p^n$ for odd primes $p\equiv 2\mod 3$. By \Cref{thm:max-sizes}, $\lambda(\F_p^n)=(p+1)p^{n-1}/3$, and it is not hard to find a sum-free set that achieves this bound. Indeed, let $m_p=(p+1)/3$ and take a subspace $V<\F_p^n$ of co-dimension 1 as well as an element $w\in \F_p^n\setminus V$. Then the set $A(V,w)=\{m_p \cdot w,\ldots, (2m_p-1)\cdot w\}+ V$ is sum-free and has size $\lambda(\F_p^n)$. Of course, any subset of $A(V,w)$ is also sum-free, and we say that a sum-free set $A$ is \emph{normal} if $A$ is contained in some set of the form $A(V,w)$. In this language, a recent result of Lev \cite{lev23} can be stated as follows.

\begin{theorem}\label{thm:p=5_lev}
If $A\subset \F_5^n$ is sum-free and $\lv A\rv > 1.5 \cdot 5^{n-1}$, then $A$ is normal.
\end{theorem}

Note that the largest possible size of a sum-free set in $\F_5^n$ is $2\cdot 5^{n-1}$. Lev expressed the view that \Cref{thm:p=5_lev} is not sharp. Our first main result confirms this.

\begin{theorem}\label{thm:F5}
If $A\subset \F_5^n$ is sum-free and $\lv A\rv > 1.2 \cdot 5^{n-1}$, then $A$ is normal.
\end{theorem}

We will briefly compare the proofs of \Cref{thm:p=5_lev} and \Cref{thm:F5} at the end of \Cref{sec:p5}. 

As part of our proof of \Cref{thm:F5}, we show that for $n=2$, the bound $6 \cdot 5^{n-2}=6$ can be replaced by $5$, which is best possible as $A=\{(0,1),(1,2),(2,2),(3,2),(4,3)\}$ shows. For larger $n$, the following example shows that the bound in \Cref{thm:F5} cannot be improved beyond $1.12\cdot 5^{n-1}$.

\begin{example}\label{example:largest-non-normal}
Consider the set 
\begin{align*}
X=\{(1,x,y)\in \F_5^3: x\in \{1,2\}\lor (x=0\land y\in \{1,2\})\}\cup \{(2,0,0), (2,0,1)\}.
\end{align*}
It is easy to see that $X$ does not intersect $X+X$, $X-X$, or $-(X+X)$, and therefore, $Y=X\cup (-X)$ is sum-free. On the other hand, $Y$ is not contained in two cosets of a 2-dimensional subspace, as such a subspace would need to contain $(0,0,1)$ and $(0,1,0)$. By taking $A=Y\times\F_5^{n-3}$ for any $n\geq 3$, we obtain a sum-free set with $28\cdot 5^{n-3}=1.12\cdot 5^{n-1}$ elements that is not normal.
\end{example}
We have no particular reason to believe that either one of the bounds $1.2\cdot 5^{n-1}$ and $1.12\cdot 5^{n-1}$ is correct. It is even conceivable that the sequence $(\lambda(\F_5^n))_{n\in \N}$ does not achieve its supremum at all. 

For larger primes $p$, we obtain the following result.

\begin{theorem}\label{thm:Fp_general}
Let $p\equiv 2\mod 3$ and $n\in \N$. If $A\subset \F_p^n$ is sum-free and $\lv A\rv\geq (m_p-1/2+1/p)p^{n-1}$, then $A$ is normal.
\end{theorem}

The bound $(m_p-1/2+1/p)p^{n-1}$ is to be compared to $\lambda(\F_p^n)=m_pp^{n-1}$. To the best of our knowledge, the previously best bound of this type is due to Green and Rusza, who proved a lemma that implies that a sum-free set $A\subset \F_p^n$ must be normal if $\lv A\rv\geq (m_p-1/(3p^2+3p))p^{n-1}$ (see Lemma 5.6. in \cite{green_rusza05}). The bound in \Cref{thm:Fp_general} cannot be sharp in so far that $(m_p-1/2+1/p)p^{n-1}$ is not an integer, and it seems likely that one can improve it by a positive proportion of $p^{n-1}$. However, the following example shows that the bound cannot be improved beyond $(m_p-1)p^{n-1}$.

\begin{example}
Let $p\geq 11$ be a prime, $p\equiv 2 \mod 3$, and let
\begin{align*}
A_2=\{(m_p-1,0),(2m_p-1,p-1)\} &\cup \{m_p\} \times (\F_p\setminus \{p-1\}) \cup \{2m_p-2\} \times (\F_p\setminus \{0\})\\ &\cup \{m_p+1,\ldots,2m_p-3\} \times \F_p\subset \F_p^2.
\end{align*}
The set $A_2\times \F_p^{n-2}\subset \F_p^n$ is sum-free, has size $(m_p-1)p^{n-1}$, and is not normal.
\end{example}

We conclude the introduction with the remark that the situation seems to be more complex for primes that are $1\mod 3$. Indeed, even for $p=7$, there are sum-free sets of size $\lambda(\F_7^n)$ that are not normal.

\begin{example}
In $\F_7^2$, the set
\begin{align*}
A=\{(3, 0), (4, 0), (3, 1), (4, 1), (3, 2), (4,2), (2,3), (3, 3), (4,3), (3, 4), (4, 4), (3, 5), (4, 5),(3,6)\}
\end{align*}
is sum-free, and has size $\lambda(\F_7^2)=14$ but is not contained in three cosets.
\end{example}

\section{Preliminaries}\label{sec:preliminiaries}
Before we state some key ingredients, let us set up some terminology. Throughout this section, let $G$ be a finite Abelian group, $p$ a prime that is $2\mod 3$ and $n\in \N$. For a set $A\subset \F_p^n$, we denote by $\langle A\rangle$ the linear span of $A$, and as usual we abbreviate $\langle \{v_1,\ldots,v_k\}\rangle$ as $\langle v_1,\ldots,v_k\rangle$ for $v_1,\ldots,v_k\in \F_p^n$. The set $I_p$ is the central interval $\{m_p+1,\ldots,2m_p-1\}\subset \F_p$, and for an element $v\in \F_p^n$ and a set $J\subset \F_p$, we write $J(v)$ for $\{xv:x\in J\}\subset \F_p^n$. 

We say that two sets $A, B\subset \F_p^n$ are \emph{isomorphic} if there exists an automorphism of $\F_p^n$, i.e., an invertible linear map $\phi\colon \F_p^n\rightarrow\F_p^n$, such that $\phi(A)=B$. Note that if $A$ is isomorphic to $B$ and $A$ is sum-free, then $B$ must be sum-free as well, and vice-versa. We will use this observation repeatedly to simplify steps of a proof by assuming that some linearly independent vectors $v_1,\ldots,v_k$ are standard unit vectors $e_1,\ldots,e_k$.

Given a set $X\subset G$, the \emph{symmetry group $\Sym(X)$ of $X$} is defined as $\{g\in G: g+X=X\}$. Note that $\Sym(X)$ is in fact a group. We will make repeated use of a theorem due to Kneser \cite{kneser53} (see \cite{tao_vu06} for a reference in English).

\begin{lemma}[Kneser's Theorem]
For a finite Abelian group $G$ and non-empty sets $A, B\subset G$, we have
\begin{align*}
\lv A + B\rv \geq \lv A+\Sym(A+B)\rv + \lv B+ \Sym(A+B)\rv - \lv \Sym(A+B)\rv.
%\geq \lv A\rv + \lv B\rv - \lv \Sym(A+B)\rv.
\end{align*}
\end{lemma}

An important special case of Kneser's Theorem is the following, which we will often use implicitly.

\begin{lemma}\label{lemma:simple_kneser}
Let $G$ be a finite Abelian group with subsets $A,B\subset G$. Suppose there exists a subgroup $H\leq G$ and $a,b\in G$ such that $A\subset a+H$  and $B\subset b+H$. If $\lv A\rv + \lv B \rv >\lv H\rv$ then $A+B=a+b+H$.
\end{lemma}
\begin{proof}
Without loss of generality, $a=b=0$, in which case clearly $A+B\subseteq H$. Suppose now that $C= H\setminus (A+B)$ is not empty. We must have that $(C-A)=C+(-A)$ is disjoint from $B$, but this is in contradiction to $\lv C-A\rv+\lv B\rv\geq \lv A\rv+\lv B\rv>\lv H\rv$.
\end{proof}

The next lemma will allow us to to gain insight into the structure of a sum-free set by considering its intersection with a subspace.

\begin{lemma}\label{lemma:subspace-finder}
Let $n\geq 3$, let $A\subset \F_p^n$ be sum-free, and let $u\in \F_p^n\setminus \{0\}$ such that $I_p(u)\subset A$. Then there exist two distinct subspaces $W_1,W_2<\F_p^n$ of co-dimension 1 such that $I_p(u)\subset W_1\cap W_2$ and $\lv W_i\cap A\rv \geq \lv A\rv /p$ for $i\in \{1,2\}$.
\end{lemma}
\begin{proof}
Let $\cW$ be the set of subspaces of $\F_p^n$ of co-dimension 1 that contain $u$. There are $K:=(p^{n-1}-1)/(p-1)$ such subspaces. We claim that least two of them intersect $A$ in more than $\lv A\rv /p$ elements. To see this, let us first determine the average size of the intersection of a subspace in $\cW$ with $A$. We have
\begin{align}\label{eq:avg-intersection1}
\sum_{W\in \cW} \lv A \cap W\rv = \sum_{w\in A} \lv \{ W\in \cW: w\in W\}\rv = Km_p + \sum_{w\in A\setminus \langle u\rangle} \lv \{W\in \cW: w\in W\}\rv.
\end{align}
All elements of $\F_p^n\setminus \langle u \rangle$ appear in the same number of subspaces in $\cW$. Since $\F_p^n\setminus \langle u \rangle$ has $p^n-p$ elements and each $W\in \cW$ contains $p^{n-1}-p$ of them, this number must be
\begin{align*}
\frac{p^{n-1}-p}{p^{n}-p}K=\frac{p^{n-2}-1}{p^{n-1}-1}K.
\end{align*}
Inserting this into \eqref{eq:avg-intersection1} yields
\begin{align}\label{eq:avg-intersection2}
\sum_{W\in \cW} \lv A \cap W\rv&=Km_p+\frac{p^{n-2}-1}{p^{n-1}-1}K(\lv A\rv -m_p)\nonumber\\
&=\frac{p^{n-1}-p^{n-2}}{p^{n-1}-1}Km_p + \frac{K\lv A \rv}{p}-\frac{p-1}{p(p^{n-1}-1)}K\lv A\rv\nonumber\\
&=m_pp^{n-2}+\frac{K\lv A\rv}{p}-\frac{\lv A\rv}{p}.
\end{align}
Assume now towards a contradiction that, but for one exception $W_1$, all subspaces in $\cW$ intersect $A$ in strictly less than $\lv A\rv/p$ points. Because $A$ is sum-free, we have $\lv A\cap W_1\rv \leq \lambda(\F_p^{n-1})= m_pp^{n-2}$, and because both $\lv A\rv$ and $\lv A\cap W\rv$ are integers for all $W\in \cW$, we have $\lv A\cap W\rv\leq (\lv A\rv-1)/p$ for all $W\in \cW\setminus \{W_1\}$. Therefore,
\begin{align*}
\sum_{W\in \cW} \lv A \cap W\rv \leq m_pp^{n-2} + \frac{(K-1)(\lv A\rv -1)}{p}=m_pp^{n-2}+\frac{K\lv A\rv}{p}-\frac{\lv A \rv}{p}+\frac{1-K}{p}.
\end{align*}
Since $K>1$ as $n\geq 3$, this is in contradiction to \eqref{eq:avg-intersection2}.
\end{proof}

Finally, the following result due to Vosper \cite{vosper56} (see also \cite{tao_vu06}) serves as the base case $n=1$ for \Cref{thm:Fp_general}.

\begin{lemma}[Vosper's Theorem]\label{lemma:vosper}
Let $p$ be a prime, and let $A, B\subset \F_p$ be such that $\lv A\rv, \lv B\rv \geq 2$ and $\lv A + B\rv  \leq p - 2$. Then $\lv A + B\rv =
\lv A\rv+\lv B\rv - 1$ if and only if $A$ and $B$ are arithmetic progressions with the same step.
\end{lemma}

\section{The proof of \Cref{thm:F5}}\label{sec:p5}

The proof of \Cref{thm:F5} involves quite a few technical steps, but in its essence, it can be split into two principal parts within an inductive framework. To give a brief overview, consider a sum-free set $A\subset \F_5^n$. Using \Cref{lemma:subspace-finder} and the induction hypothesis, we will find a subspace $V<\F_5^n$ such that $A\cap V$ is normal, i.e., contained in two cosets of some $(n-2)$-dimensional subspace $U<V$. 

In the first part of the induction step, we will show that if $A$ is not sufficiently uniformly distributed over the other cosets of $U$ in $\F_5^n$, then $A$ can be reduced to a set $A'\subset \F_5^2$ that is also sum-free. Building on a good understanding of the 2-dimensional case, this will allow us to infer structural properties of $A'$ from which we will be able to deduce that $A$ must be normal. In the second part of the induction step, we will exploit the fact that $A$ is sum-free in a series of linear inequalities to show that $A$ must in fact display a highly irregular distribution over the cosets of $U$, thus completing the proof.

Although the separation between the two parts is not entirely clean, the first part roughly maps to \Cref{lemma:2-dim-shape}, \Cref{lemma:rich-triple}, \Cref{claim:fat-set-six}, and \Cref{claim:combined-five}, while the remainder of the proof of \Cref{thm:F5} covers the second part. 

\begin{lemma}\label{lemma:2-dim-shape}
%If a sum-free set $A\subset \F_5^2$ of size 5 is not normal, then $A$ is isomorphic to either $\{(0,1),(1,1),(2,1),(2,3),(0,4)\}$ or $\{(0,1),(1,1),(2,1),(4,3),(3,4)\}$.
Let $A\subset \F_5^2$ be a sum-free set of size at least 5. Then $A$ is isomorphic to a set that contains $(1,0)$, $(1,1)$ and $(1,2)$.
\end{lemma}
\begin{proof}
It is enough to show that there exists a 1-dimensional subspace $W< \F_5^2$ and $u\in \F_5^2$ such that $\lv A\cap (u+W)\rv\geq 3$. Suppose this is not the case. Since $A$ is sum-free, we must have $A\cap (A-A)=\emptyset$. Therefore, we must have $\lv A-A\rv \leq 20$, which means that there must be $a,b,c,d\in A$ such that $a\neq b,c$ and $a-b=c-d$. We cannot have $a=d$ or $b=c$ either since then there would be three points in arithmetic progression. It remains to check that there is no $f\in \F_5^2\setminus \{a,b,c,d\}$ such that $\{a,b,c,d,f\}$ is sum-free and no two points of $\{a,b,c,d\}$ are collinear with $f$. To simplify this process, we may assume that $a-b=c-d=e_1$ and $a-c=e_2$.  Writing $d$ as $(x,y)$ for $x,y\in \F_5$, we have $c=(x+1,y), a=(x+1,y+1)$ and $b=(x,y+1)$. By symmetry, we may assume that $x\geq y$ when viewed as elements of $\{0,\ldots,4\}$ and because $A$ is sum-free, we must have that $(x,y)\notin \{(0,0),(1,0),(3,0),(4,0),(1,1),(4,1),(3,3),(4,3),(4,4)\}$. This leaves six possible choices for $(x,y)$ to check by hand, a task we leave to the reader.
\end{proof}

The following two lemmas are direct consequences of the above.

\begin{lemma}\label{lemma:2-dim-empty-cosets}
If a sum-free set $A\subset \F_5^2$ has size at least 5, then there exists a 1-dimensional subspace $V< \F_5^2$ such that $(A+A)\cup (A-A)$ contains two cosets of $V$.
\end{lemma}
\begin{proof}
By \Cref{lemma:2-dim-shape}, we may assume without loss of generality that $(1,0),(1,1),(1,2)\in A$, whence $\langle e_2\rangle \subset A-A$ and $2e_1+\langle e_2\rangle \subset A+A$.
\end{proof}

\begin{lemma}\label{lemma:F5-base-case}
\Cref{thm:F5} holds for $n=1$ and $n=2$.
\end{lemma}
\begin{proof}
The claim is trivial for $n=1$, and for $n=2$, it is readily deduced from \Cref{lemma:2-dim-empty-cosets}.
\end{proof}

The next lemma applies \Cref{lemma:subspace-finder} to the case $p=5$ and allows us to prove \Cref{thm:F5} by induction.

\begin{lemma}\label{lemma:subspace-finder-F5}
For $n\geq 2$ and any sum-free set $A\subset \F_5^n$ with $\lv A\rv > 1.2\cdot 5^{n-1}$, there is a subspace $W< \F_5^n$ of co-dimension 1 such that $\lv A\cap W\rv \geq \lv A \rv/5$.
\end{lemma}
\begin{proof}
Suppose first that $n=2$. By assumption, $\lv A\rv\geq 6$, and so by \Cref{lemma:F5-base-case}, $A$ must be normal. Without loss of generality $A\subset \{e_1,-e_1\}+\langle e_2\rangle$, and there must exist $y\in \F_5$ such that $(1,y),(-1,-y)\in A$.

Suppose now that $n\geq 3$. By \Cref{lemma:subspace-finder}, it is sufficient to find $u\neq 0$ such that $I_p(u)\subset A$. Using the argument above, it is enough to find a 2-dimensional subspace $V<\F_5^n$ such that $\lv A\cap V\rv\geq 6$. For this purpose, let $K$ be the number of 2-dimensional subspaces of $\F_5^n$. Every element of $A$ appears in $K(5^2-1)/(5^n-1)$ of them. Therefore, a 2-dimensional subspace contains on average
\begin{align*}
\frac{\lv A\rv(5^2-1)}{5^n-1}>\frac{28.8}{5}
\end{align*}
elements of $A$, which means one of them contains at least six, as desired.
\end{proof}

We remark that with a little more effort, one can improve the bound $1.2\cdot 5^{n-1}$ in \Cref{lemma:subspace-finder-F5} to $5^n/6$.

\begin{lemma}\label{lemma:rich-triple}
Let $n\geq 2$ and $A\subset \F_5^n$ be sum-free. If there exists a subspace $W<\F_5^n$ of co-dimension 1, $v\in \F_5^n$ and $x_1,x_2,x_3\in \F_5$ such that $\lv A\cap (\{x_1v,x_2v,x_3v\}+W)\rv > 26 \cdot 5^{n-3}$, then $A$ is normal.
\end{lemma}
\begin{proof}
Without loss of generality, $W$ is the coordinate plane $\{w: w_1=0\}$ and $v=e_1$. For $x\in \F_5$, we let $C_x= xe_1+W$ and $A_x=A\cap C_x$. Suppose first that one of $A_{x_1}, A_{x_2}$ and $A_{x_3}$,  $A_{x_3}$ say, is empty. In this case, we have $\lv A_{x_1} \rv + \lv A_{x_2}\rv > 5^{n-1}$, whence we must have $A_{x_1}-A_{x_2}=C_{x_1-x_2}$, $A_{x_2}-A_{x_1}=C_{x_2-x_1}$ and $(A_{x_1}-A_{x_1})\cup (A_{x_2}-A_{x_2})=C_{0}$. Hence, $A$ is normal.

Assume now that none of $A_{x_1}, A_{x_2}$ and $A_{x_3}$ are empty, and without loss of generality, that $x_1+x_2=x_3$. We have that $\lv A_{x_1}\rv+\lv A_{x_2}\rv + \lv A_{x_3}\rv > 26 \cdot 5^{n-3}$, but for $A$ to be sum-free, we must have $\lv A_{x_1}+A_{x_2}\rv + \lv A_{x_3}\rv \leq 5^{n-1}$. By Kneser's Theorem, this is only possible if the symmetry group $U:= \Sym(A_{x_1} + A_{x_2})$ has size $5^{n-2}$. Since $ A_{x_1} + A_{x_2} \subset C_{x_3}$, we must have $U\subset W$.

Again without loss of generality, $e_2\notin U$, and we consider the set $S=\{(x,y)\in F_5^2: x\in \{x_1,x_2,x_3\} \text{ and } A\cap (U+xe_1+ye_2)\neq \emptyset\}$. We know that $\lv S \rv$ must be at least 6, since five cosets of $U$ could fit at most $25\cdot 5^{n-3}$ elements. On the other hand, no coset of $U$ can contain an element of both $A_{x_1} + A_{x_2}$ and $A_{x_3}$, and therefore it is easy to see that $\lv S\rv=6$. 

Suppose that $S$ is normal. In this case, there must exists $a\in \F_5^2$ such that $S$ is contained in two cosets of $\langle a\rangle$. Letting $V=\langle U, a_1e_1+a_2e_2\rangle$, we see that $A_{x_1}\cup A_{x_2} \cup A_{x_3}$, which contains more than $26\cdot 5^{n-3}$ elements, is contained in two cosets of $V$. Repeating the first step of the proof, we see that $A$ must be normal.

If $S$ is not normal, \Cref{lemma:F5-base-case} tells us that $S$ cannot be sum-free. Let $a,b,c\in S$ such that $a+b=c$. Since we may assume that at most two of the three points $a,b$ and $c$ can be equal, and the union of all cosets of $U$ with coordinates in $S$ intersects $A$ in at least $26\cdot 5^{n-3}$ elements, we must have $\sum_{(x,y)\in \{a,b,c\}}\lv  A\cap (U+xe_1+ye_2)\rv>7\cdot 5^{n-3}$, whence $A$ cannot be sum-free.
\end{proof}

We are now ready to prove \Cref{thm:F5} in full.

\begin{proof}[Proof of \Cref{thm:F5}]
We prove the claim by induction on $n$. Assume therefore that the claim is true for $n-1$, and let $A\subset \F_5^n$ be a sum-free set that satisfies $\lv A\rv > 6\cdot 5^{n-2}$ but is not normal.

For a subspace $V\leq \F_5^n$ and an element $v\in \F_5^n$, we define $d_{V,v}$ as $5\lv A \cap (v+V)\rv/\lv V\rv$. We claim that there is an $(n-2)$-dimensional subspace $U$ and an element $u\in \F_5^n$ such that $d_{U,u}>3$. Indeed, note that by \Cref{lemma:subspace-finder-F5}, there is a subspace $V< \F_5^n$ of co-dimension 1 such that $d_{V,0}> 1.2$, and since $A\cap V$ is sum-free, we can conclude from the induction hypothesis that $V$ must have a subspace $U\leq V$ of dimension $n-2$ such that $A\cap V$ is contained in two cosets of $U$. At least one of these cosets must intersect $A$ in more than $3\cdot 5^{n-3}$ elements, as desired.

We now choose $U$ of dimension $n-2$ and $v \in \F_5^n$ such that $d_{U,v}$ is maximal. Since $A$ is sum-free, and $d_{U,v}>3$, we must have $v\neq 0$, so that $V=\langle U\cup \{v\}\rangle$ has dimension $n-1$. Without loss of generality, $v=e_1$ and $e_2\notin V$. In what follows, we will think of cosets of $U$ as being arranged in a 2-dimensional grid, referring to them as \emph{cells} and to cosets of $V$ as \emph{rows}. For $x,y\in \F_5$, we define $U_{x,y}:=U+xe_1+ye_2$, $A_{x,y}=A\cap U_{x,y}$ and $f(x,y)=\lv A_{x,y}\rv/5^{n-3}$. For a set $S\subset \F_5^2$, we define $f(S)=\sum_{(x,y)\in S} f(x,y)$.

\begin{claim}\label{claim:fat-set-six}
The set $S=\{(x,y):f(x,y)>2\}$ has size at most 5. 
\end{claim}

\begin{proof}[Proof of \Cref{claim:fat-set-six}.]
Suppose towards a contradiction that $\lv S\rv \geq 6$. By Kneser's Theorem, $S$ must be sum-free, and thus by \Cref{lemma:F5-base-case}, $S$ must be normal. We may assume without loss of generality that $S\subset (e_1+\langle e_2\rangle)\cup (-e_1+\langle e_2\rangle)$ while maintaining that $f(1,0)>3$. Let $T=\{(x,y)\notin (e_1+\langle e_2\rangle)\cup (-e_1+\langle e_2\rangle): f(x,y)>0\}$. Because $f(1,0)>3$, we must have $f(x+1,y)=f(x-1,y)=0$ for all $(x,y)\in S$, and thus it is easy to see that $\lv T\rv$ is at most 6. 
We also know that $T\subset (S+S)\cup (S-S)$, and for all $(x_1,y_1),(x_2,y_2)\in S$ we have $\vert A_{x_1,y_1} \pm A_{x_2,y_2}\vert > 4 \cdot 5^{n-3}$ by Kneser's Theorem, whence $f(x,y)<1$ must hold for all $(x,y)\in T$.

Since $f(\F_5^2)>30$ and $f(T)<6$,  there must exist $x\in \{-2,0,2\}$ such that $f(\{e_1,-e_1,xe_1\}+W)>26$, which is excluded by \Cref{lemma:rich-triple}.
\end{proof}

\begin{claim}\label{claim:combined-five}
There is no set $\lv S\rv\subset \F_5^2$ of size at least 5 such that $f(x,y)>1$ for all $(x,y)\in S$ and $f(x_1,y_1)+f(x_2,y_2)>5$ for all distinct $(x_1,y_1),(x_2,y_2)\in S$.
\end{claim}
\begin{proof}[Proof of \Cref{claim:combined-five}.]
Assume towards a contradiction that a set $S$ such as in the claim exists. By Kneser's Theorem, $S$ must be sum-free, and by \Cref{lemma:2-dim-shape}, we may therefore assume without loss of generality that $S$ contains the points $(1,0),(1,1)$ and $(1,2)$. Let $y_0\in \{0,1,2\}$ such that $f(1,y_0)$ is minimal. By the second assumption on $S$ in the claim, we must have $f(1,y)> 2.5$ for $y\in \{0,1,2\}\setminus \{y_0\}$. By Kneser's Theorem, we have that $f(0,y)=0$ for all $y\in \F_5$, $f(2,y)=0$ for all $y\in \F_5\setminus \{2y_0\}$, and $f(2,2y_0)<3$. Thus, $f(\{e_1,3e_1,4e_1\}+\langle e_2\rangle)>27$, contradicting \Cref{lemma:rich-triple}.
\end{proof}

We say that a row $ye_2+V$ is \emph{scattered} if there are at least four cells contained in $ye_2+V$ that have non-empty intersection with $A$. Similarly, we say that a row is \emph{semi-scattered} or \emph{focused} if the number of such cells is exactly three or at most two, respectively.

\begin{claim}\label{claim:scattered-coset4}
Suppose $f(1,0)>4$. Every scattered row $ye_2+V$ satisfies $f(ye_2+\langle e_1\rangle)\leq 2.5$, and if $ye_2+V$ is semi-scattered, then we can find $(x_1,y),(x_2,y)$ such that $f(x_1,y), f(x_2,y)>0$ and $f(x_1,y)+f(x_2,y)\leq 1$.
\end{claim}
\begin{proof}[Proof of \Cref{claim:scattered-coset4}.]
Suppose there is a scattered row $ye_2+V$ such that $f(ye_2+\langle e_1\rangle)>2.5$. By an averaging argument, there must exist $x\in \F_5$ such that $f(x,y)+f(x+1,y)>1$. One of $A_{x,y}$ and $A_{x+1,y}$ must be empty, as otherwise by Kneser's Theorem,
\begin{align*}
\lv A_{x,y} + A_{1,0}\rv + \lv A_{x+1,y} \rv> \lv A_{x,y}\rv + \lv A_{x+1,y}\rv + 4\cdot 5^{n-3}>5^{n-2},
\end{align*}
which would contradict the assumption that $A$ is sum-free. But then we can assume without loss of generality that $f(x,y)>1$ and conclude by the same argument that $A_{x-1,y}$ must be empty as well, thus contradicting the assumption that $ye_2+V$ is scattered.

To show the second part of the claim, we simply find two among the three cells with non-empty intersection with $A$ in $ye_2+V$ whose $x$-coordinates differs by 1 and apply Kneser's Theorem to them.
\end{proof}

\begin{claim}\label{claim:scattered-coset3}
Every scattered row $ye_2+V$ satisfies $f(ye_2+\langle e_1\rangle)\leq 5$, and if $ye_2+V$ is semi-scattered, we can find $(x_1,y),(x_2,y)$ such that $f(x_1,y), f(x_2,y)>0$ and $f(x_1,y)+f(x_2,y)\leq 2$.
\end{claim}
\begin{proof}[Proof of \Cref{claim:scattered-coset3}.]
Same as \Cref{claim:scattered-coset4}, noting that $f(1,0)>3$.
\end{proof}

\begin{claim}\label{claim:two-scattered-rows}
There can be at most two scattered rows.
\end{claim}
\begin{proof}[Proof of \Cref{claim:two-scattered-rows}.]
A single cell can contain at most $5\cdot 5^{n-2}$ elements of $A$, whence a focused row can contain at most $10\cdot 5^{n-2}$ elements of $A$, and by \Cref{claim:scattered-coset3}, a semi-scattered row can contain at most $7\cdot 5^{n-2}$ elements of $A$. Since every scattered row intersects $A$ in at most $5\cdot 5^{n-3}$ elements and $\lv A\rv > 30 \cdot 5^{n-3}$, we know that there can be at most three scattered rows, and if there are indeed three of them, then at least one of the remaining ones has to be focused. Furthermore, $V$ can never be scattered as $f(0,0)=f(2,0)=0$. Suppose now towards a contradiction that there are only two non-scattered rows, namely $V$ and $ye_2+V$ for some $y\in \F_5$.

We assume first that $V$ and $ye_2+V$ are both focused. In $V\cup (ye_2+V)$, there are at most four cells that have a non-empty intersection with $A$, and since $\lv A\setminus (V\cup (ye_2+V))\rv\leq 15 \cdot 5^{n-3}$, this means that at least three cells within $V\cup (ye_2+V)$ have an intersection with $A$ that is larger than $2.5\cdot 5^{n-3}$. We fix coordinates $x_1,x_2,x_3,y_1$ and $y_2$ such that $x_1\neq x_2$, $y_1\neq y_2$, and $\lv A_{x_1,y_1} \rv, \lv A_{x_2,y_1} \rv, \lv A_{x_3,y_2} \rv\geq 2.5\cdot 5^{n-3}$. Without loss of generality $y_1-y_2\notin \{y_1,y_2\}$, and thus $A_{y_1-y_2,x_1-x_3}=A_{y_1-y_2,x_2-x_3}=\emptyset$, which means that $(y_1-y_2)e_2+V$ is not scattered.

If one of $V$ or $ye_2+V$ is semi-scattered, without loss of generality the latter, there exists by \Cref{claim:scattered-coset3} an element $x\in \F_5$ such that $\lv A\cap (ye_2+V)\setminus A_{x,y}\rv \leq 2\cdot 5^{n-3}$. This means that $\lv A_{x,y} \cup (A\cap V)\rv > 13\cdot 5^{n-3}$ so that there must be three cells in $V\cup (yw+V)$ that intersect $A$ in at least $3\cdot 5^{n-3}$ elements, and we can proceed as before.
\end{proof}

We are now going to analyze $f$ by partitioning $\F_5^2$ into suitable sets. Firstly, we denote by $R$ the coordinate pairs of cells that lie within a scattered row. Secondly, we take a set $S\subset \F_5^2$ that contains for each semi-scattered row $ye_2+V$ two pairs of coordinates $(x_1,y)$ and $(x_2,y)$ such that $f(x_1,y), f(x_2,y)>0$ and $f(x_1,y)+f(x_2,y)\leq 1$ if $f(1,0)> 4$ or $f(x_1,y)+f(x_2,y)\leq 2$ if $f(1,0)\leq 4$. Such a set exists by \Cref{claim:scattered-coset4} and \Cref{claim:scattered-coset3}. We denote by $r$ and $s$ the number of scattered and semi-scattered rows, respectively. 

Now, let $T=\{(x,y)\in \F_5^2\setminus (R\cup S): f(x,y)>0\}$. Each focused line can contain at most two cells with coordinates in $T$, and therefore we have $\lv T\rv \leq 10-2r-s$. Suppose there are $\e \geq 0$ and a set $E\subset T$ such that $f(E)\leq \e$, and for $c>0$, let $K_c=\{(x,y)\in T:f(x,y)>c\}$. For all $c>0$, we have
\begin{align*}
f(\F_5^2)&=f(R)+f(S)+f(E)+f(K_c\setminus E)+f(T\setminus K_c\setminus E)
\end{align*}
Recalling that $f(1,0)$ maximizes $f$ and $f(\F_5^2)> 30$, we obtain 
\begin{align}\label{eq:mass-bound-general}
 f(R)+f(S) + \e + f(1,0)\lv K_c\setminus E\rv + c(10-2r-s-\lv E\rv-\lv K_c\setminus E\rv)> 30.
\end{align}
\begin{claim}
We must have $f(1,0)>4$.
\end{claim}
\begin{proof}
Observe first that by \Cref{claim:scattered-coset3}, $f(R)\leq 5r$ and $f(S)\leq 2s$. Suppose towards a contradiction that $f(1,0)\leq 4$. Inserting this into \eqref{eq:mass-bound-general} for $c=2$, $E=\emptyset$ and $\e=0$, we obtain
\begin{align*}
r + 2\lv K_2\rv + 20 >30.
\end{align*}
Since $r$ is at most 2 by \Cref{claim:two-scattered-rows}, this implies that $\lv K_2\rv$ is more than 4, hence at least 5. By \Cref{claim:combined-five}, there must be two points $(x_1,y_1),(x_2,y_2)\in K_2$ such that $f(x_1,y_1)+f(x_2,y_2)\leq 5$. Taking $E=\{(x_1,y_1), (x_2,y_2)\}$ and $\e=5$, we can apply \eqref{eq:mass-bound-general} again with $c=2$ to obtain
\begin{align*}
r + 5 + 2\lv K_2\setminus E\rv + 16 >30,
\end{align*}
from which it follows that $\lv K_2\setminus E\rv$ is at least 4. We arrive at a contradiction by applying \Cref{claim:fat-set-six} to $K_2$.
\end{proof}

To conclude the proof of \Cref{thm:F5}, let us now analyze what happens if $f(1,0)>4$. We know by \Cref{claim:scattered-coset4} that $f(R)\leq 2.5r$ and $f(S)\leq s$. Suppose first that there exists no $(x,y)\in T$ such that $f(x,y)\in [1.5,3]$. Then we may apply \eqref{eq:mass-bound-general} with $c=1.5$, $E=\emptyset$ and $\e=0$ to obtain
\begin{align*}
3.5\lv K_{1.5}\rv + 15>30,
\end{align*}
which implies that $\lv K_{1.5}\rv\geq 5$. However, since $K_{1.5}=K_3$, we arrive at a contradiction via \Cref{claim:combined-five}.

Suppose next that there exists $(x,y)\in T$ such that $f(x,y)\in (2,3]$. We may apply \eqref{eq:mass-bound-general} with $E=\{(x,y)\}, \e=3$ and $c=2$ to obtain
\begin{align*}
3+3\lv K_2\setminus E\rv+18>30,
\end{align*}
which implies that $\lv K_2 \setminus E\rv \geq 4$. By \Cref{claim:combined-five}, there must be $(x_1,y_1),(x_2,y_2)\in K_2$ such that $f(x_1,y_1)+f(x_2,y_2)\leq 5$. Applying \eqref{eq:mass-bound-general} once more with $E=\{(x_1,y_1),(x_2,y_2)\}$, $\e=5$ and $c=2$, we obtain
\begin{align}\label{eq:exception-set-5}
5+3\lv K_2\setminus E\rv + 16 >30,
\end{align}
which implies that $\lv K_2\setminus E\rv\geq 4$. We arrive at a contradiction by applying \Cref{claim:fat-set-six} to $K_2$.

This leaves the possibility that $f$ takes values in $[1.5,2]$ but not in $(2,3]$. We define
\begin{align*}
\alpha = \max \{f(x,y):(x,y)\in T\setminus K_2\} \qquad \text{ and } \qquad \beta = \min \{f(x,y):(x,y)\in K_2\}.
\end{align*}
If $\alpha+\beta\leq 5$, we fix $(x_1,y_1)$ and $(x_2,y_2)$ such that $f(x_1,y_1)=\alpha$ and $f(x_2,y_2)=\beta$. Taking $E=\{(x_1,y_1),(x_2,y_2)\}$, we obtain again \eqref{eq:exception-set-5} so that $\lv K_2\setminus E\rv \geq 4$, which allows us to apply \Cref{claim:combined-five} to $K_2=K_3$.

We may there assume that $\alpha+\beta > 5$. Applying \eqref{eq:mass-bound-general} one last time with $E=\emptyset$, $\e=0$ and $c=2$, we obtain
\begin{align*}
3\lv K_2\setminus E\rv + 20 >30,
\end{align*}
which implies that $\lv K_2\rv\geq 4$. Here, we apply \Cref{claim:combined-five} to $K_2 \cup \{(x_1,y_1)\}$ which yields again a contradiction.
\end{proof}

Before we move on to larger primes, we briefly discuss how the proof of \Cref{thm:F5} relates to Lev's proof of \Cref{thm:p=5_lev} in \cite{lev23}. Lev's proof can be split into the same two main parts that we have laid out at the start of the section, i.e., one part in which one finds a large subspace so that $A$ is distributed non-uniformly over its cosets and one part that exploits this property to show that $A$ is normal.

The key difference between Lev's proof and ours is that in \cite{lev23}, said large subspace is $(n-1)$-dimensional and obtained through a Fourier-theoretic argument, whereas we use an inductive argument to find an $(n-2)$-dimensional subspace. The latter approach leads to improvements in both parts of the proof. Firstly, the inductive argument yields strong structural information about $A$ even when $A$ is small. In contrast to that, the quantitative information that can be extracted from Fourier coefficients quickly diminishes as the size of $A$ decreases. Secondly, considering an $(n-2)$-dimensional subspace allows us to exploit insights such as \Cref{lemma:2-dim-shape} about the structure of 2-dimensional sum-free sets. 

It seems likely that \Cref{thm:F5} can be improved further by analyzing the structure of non-normal sum-free sets in $\F_5^k$ for fixed small $k\geq 3$. However, the complexity of this analysis increases rapidly in $k$, and it is unclear whether this approach would ever yield an exact bound for the size of a largest non-normal sum-free subset of $\F_5^n$ for all $n$.

\section{The proof of \Cref{thm:Fp_general}}\label{sec:p_general}
For the entirety of this section, let $p\geq 11$ be a prime that is $2\mod 3$, $I=I_p$, $m=m_p$ and $\alpha = 1/2-1/p$. Instead of proving \Cref{thm:Fp_general} directly, we will prove the following proposition, which is more convenient for our proof by induction.

\begin{proposition}\label{prop:Fp_general}
Let $A\subset \F_p^n$ be sum-free and $u\in A$ such that $I(u)\subset A$ and $\lv A\rv\geq (m-\alpha)p^{n-1}$. There exists a subspace $V$ of co-dimension 1 such that $A\subset I(u) + V$.
\end{proposition}

\begin{proof}[Proof of \Cref{thm:Fp_general} from \Cref{prop:Fp_general}].
Let $n\in \N$ and let $A\subset \F_p^n$ be sum-free and have size at least $ (m-1/2+1/p)p^{n-1}$. Clearly, it is enough to find $u\in \F_p^n$ such that $I(u)\subset A$. There are $(p^n-1)/(p-1)$ 1-dimensional subspaces of $\F_p^n$, and we observe that every $v\in A$ is contained in exactly one of them. The average 1-dimensional subspace of $\F_p^n$ therefore contains
\begin{align*}
\frac{(p-1)\lv A \rv}{p^n-1}>\frac{(p-1)(p-1/2+1/p)}{3p}>\frac{p-2}{3}
\end{align*}
elements of $A$. We can thus find $v\in \F_p^n$ such that $A'= A\cap \langle v\rangle$ satisfies $\lv A'\rv = m$. Since $(A'+A')\cap A'=\emptyset$, we have $\lv A'+A'\rv\leq 2m-1$. On the other hand, $\lv A'+A'\rv \geq 2m-1$ by Kneser's Theorem, and so by Vosper's Theorem (\Cref{lemma:vosper}), $A'$ must be in arithmetic progression in $\langle v\rangle$. By scaling $v$ by a suitable factor, we get $u$ such that $I(u)=A'\subset A$.
\end{proof}

To prove \Cref{prop:Fp_general}, we require the following elementary lemma.

\begin{lemma}\label{lemma:continuity}
Let $d\in \{0,\ldots,m-1\}$ and $J\subset \{m-d,\ldots,2m-1+d\}$ be a subset of $\F_p$ that has size at least $d+1$ and is contained in a translate of $I$. Then $\F_p\setminus (I+J)\subset \{\lv J\rv - d,\ldots, d-\lv J\rv\}$.
% Used to be $\{2m+d+(m-\lv J\rv),\ldots,0,\ldots, m-d-(m-\lv J\rv) -1\}\subset I+J$.
\end{lemma}
\begin{proof}
Since $J$ is contained in a translate of $I$, $I+J$ is an interval. Let $x$ be the minimal element in $J$ considered as a subset of $\{0,\ldots,p-1\}$. We know that $x\geq m-d>0$ and $x\leq 2m+d-\lv J\rv\leq 2m-1$. Therefore, $I+J$ contains
\begin{align*}
\{2m+d-\lv J\rv+m,\ldots,0,\ldots,m-d+\lv J \rv -1+2m-1\}=\{ d-\lv J\rv + 1,\ldots, \lv J\rv - d -1\},
\end{align*}
and the lemma follows by taking complements.
\end{proof}

\begin{proof}[Proof of \Cref{prop:Fp_general}]
We prove the claim by induction on $n$, where the base case $n=1$ is trivial. Suppose therefore that $n\geq 2$, that $A\subset\F_p^n$ is sum-free, and, without loss of generality, that $I(e_1)\subset A$. For an element $v\in \F_p^n$, we define $A_v$ to be the unique subset of $\F_p$ such that $A_v(e_1)=(A-v)\cap\langle e_1\rangle$. In particular, $A_0=I$ and $(A_0+A_v)\cap A_v=\emptyset$ for all $v\in \F_p^n$. By Kneser's Theorem, this implies that $\lv A_v\rv\leq m$ for all $v$, but at the same time, we have that the average size of $A_v$ is at least $m-\alpha$. This will allow us to conclude that $\lv A_v\rv \geq m$ for most $v$. We will call such $v$ \emph{good}, and call $v$ \emph{bad} otherwise. 

Observe that if $v$ is good, then by Kneser's Theorem and Vosper's Theorem, $A_v$ must be a translate of the interval $I$, i.e., $A_v=t_v+I$ for some $t_v\in \F_p$. What is more, there are certain relations between these offsets $t_v$ for different $v$, which are captured by the following claim.
\begin{claim}
For all good $v,w$ we have 
\begin{align}\label{eq:offset-addition}
A_{v+w}&\subset t_v+t_w+I
\end{align}
and
\begin{align}\label{eq:offset-subtraction}
A_{v-w}\subset t_v-t_w+I.
\end{align}
\end{claim}
\begin{proof}
Noting that $I+I=I-I=\F_p\setminus I$, we have
\begin{align*}
A_{v+w}\subset \F_p \setminus (A_v+A_w)=t_v+t_w+( \F_p\setminus (I+I))=t_v+t_w+I
\end{align*}
and
\begin{align*}
A_{v-w}\subset \F_p \setminus (A_v-A_w)=t_v-t_w+( \F_p\setminus (I-I))=t_v-t_w+I.
\end{align*}
\end{proof}

Since there are many good $v$, we will be able to conclude that even for bad $v$, the set $A_v$ is \emph{contained} in some translate of $I$. The remainder of the proof is devoted to making sure that these translates actually line up relative to some subspace $V$ of co-dimension 1. As a first step to this end, we generate a suitable candidate for $V$, and our approach to this task depends on whether $n$ is 2 or larger than 2. We begin with the former case.

\begin{claim}\label{claim:avoided-subspace-n2}
Suppose $n=2$. There exists a subspace $V\subset \F_p^2$ of dimension 1 such that $A\cap V=A\cap (e_1+V)=A\cap(-e_1+V)=\emptyset$.
\end{claim}
\begin{proof}[Proof of \Cref{claim:avoided-subspace-n2}.]
Take an arbitrary element $v\in \F_p^2\setminus \langle e_1\rangle$ and consider the set $G=\{x\in \F_p: \lv A_{xv}\rv =m\}$. By a simple counting argument, the size of $G$ is at least $p/2+1$, so there must exist a non-zero $x$ such that both $A_{xv}$ and $A_{2xv}$ are in $G$. We may assume without loss of generality that $x = 1$ and that $t_v = 0$, which can always be achieved by first scaling $v$ and then translating it by a suitable multiple of $e_1$. By \eqref{eq:offset-addition}, this also implies that $t_{2v}=2t_v=0$.

We partition $\F_p^*$ into $H=\{2,4,\ldots,p-1\}$ and $H-1$. Suppose that $\sum_{x\in H} \lv A_{xv}\rv\geq \sum_{x\in (H-1)} \lv A_{xv}\rv$. In particular, we must have
\begin{align*}
\sum_{x\in H} m-\lv A_{xv}\rv\leq \frac{(p-1)m}{2}-\frac{(m-\alpha)p-m}{2} = \frac{\alpha p}{2}< m-1,
\end{align*}
and since the sum on the left hand side is an integer, it must be bounded by $m-2$. Considering $H$ as a subset of $\N$, we define
\begin{align*}
f\colon H\rightarrow \N_0 \qquad x \mapsto 
\sum_{\substack{z \in H\\z <x}} m-\lv A_{zv} \rv.
\end{align*}
Note that for all $x\in H$ we have $f(x)\leq m-2-(m-\lv A_{xv})\rv=\lv A_{xv}\rv-2$ and that since $\lv G\rv>p/2$, $A_{xv}$ must be contained in a translate of $I$ for all $x\in \F_p$ by \eqref{eq:offset-addition} and \Cref{lemma:simple_kneser}. Using \Cref{lemma:continuity} and the fact that $A_{2v}=I$, we can thus prove inductively that $A_{xv} \subset \{m-f(x),\ldots,2m-1+f(x)\}$ for all $x \in H$. Finally, using that $A_v=I$ as well, we can apply \Cref{lemma:continuity} once more to $A_{xv}$ for all $x\in H$ to see that $A_{x-1}\subset \{2,\ldots,3m-3\}$. Taking $V=\langle v\rangle$ thus gives a subspace as desired. If $\sum_{x\in (H-1)} \lv A_{xv}\rv> \sum_{x\in H} \lv A_{xv}\rv$, we may switch the roles of $H$ and $H-1$ and argue analogously.
\end{proof}
For $n\geq 3$, the argument is slightly longer, but most of it comes down to rather simple linear algebra.
\begin{claim}\label{claim:avoided-subspace-n3}
Suppose $n\geq 3$. There exists a subspace $V\subset \F_p^n$ of co-dimension 1 such that $A\cap V=A\cap (e_1+V)=A\cap(-e_1+V)=\emptyset$.
\end{claim}
\begin{proof}[Proof of \Cref{claim:avoided-subspace-n3}.]
By \Cref{lemma:subspace-finder}, there are two distinct subspaces $W_1$ and $W_2$ of co-dimension 1 such that $\langle e_1\rangle \subset W_1\cap W_2$ and $\lv A\cap W_1\rv, \lv A\cap W_2\rv \geq \lv A\rv /p$. For $i\in \{1,2\}$, we can apply the induction hypothesis to $A\cap W_i$ to see that there exists an $(n-2)$-dimensional subspace $V_i<W_i$ such that $A\cap W_i$ is contained in $I(e_1)+V_i$. Note that $\lv A\cap W_i\rv$ is larger than $\lv V_i\rv$ so that $e_1$ cannot be contained in $V_i$, and since $e_1\in W_1\cap W_2$, we can conclude that $\dim V_1\neq V_2$.

Suppose now that $e_1\in V_1+V_2$. If follows that $W_1, W_2\subset V_1+V_2$, and therefore, that $(x_1e_1+V_1)+(x_2e_1+V_2)=\F_p^n$ for all $x_1,x_2\in \F_p$. By an averaging argument, there must exist $x_1,x_2\in \F_p$ such that for $i\in \{1,2\}$,
\begin{align*}
\lv A\cap (x_ie_1+V_i)\rv \geq \frac{\lv A \rv}{mp} \geq p^{n-2}(1-\alpha/m).
\end{align*}
At the same time, $V_1+V_2=\F_p^n$ also implies that $\dim V_1\cap V_2=n-4$ so that
\begin{align*}
\lv A\cap (x_1e_1+V_1)\rv\lv A\cap (x_2e_1+V_2)\rv \leq \lv V_1\cap V_2\rv \lv A+A\rv = p^{n-4}\lv A+A\rv.
\end{align*}
We thus obtain that $\lv A+A\rv\geq (1-\alpha/m)^2p^n>3p^n/4$, which cannot hold since $A$ is sum-free and $\lv A\rv >p^n/4$. Consequently, we must have $e_1\notin V_1+V_2$ and therefore $\dim V_1\cap V_2=n-3$. Letting $V=V_1+V_2$, we conclude that $\dim V=n-1$.

For $v\in V$, let $v^{(1)}\in V_1,v^{(2)}\in V_2$ be such that $\lv A_{v^{(1)}}\rv +\lv A_{v^{(2)}}\rv$ is maximal with $v^{(1)}+v^{(2)}=v$. Certainly, the maximum must at least as large as the average, so that
\begin{align}\label{eq:v12-against-avearge}
\lv A_{v^{(1)}}\rv +\lv A_{v^{(2)}}\rv \geq p^{-(n-3)}\sum_{\substack{w^{(1)}\in V_1, w^{(2)}\in V_2\\w^{(1)}+w^{(2)}=v}} \lv A_{w^{(1)}} \rv + \lv A_{w^{(2)}}\rv.
\end{align}
On the other hand, we know that $A_v\subset  \F_p\setminus (A_{v^{(1)}}+A_{v^{(2)}})$ and that $A_w\subset I$ for all $w\in V_1\cup V_2$, which yields (by Kneser's Theorem for $\F_p$)
\begin{align*}
\lv A_v\cap I\rv\geq \lv A_v\rv - (2m - 1 - \lv A_{v^{(1)}}+A_{v^{(2)}}\rv)\geq \lv A_v\rv + \lv A_{v^{(1)}} \rv + \lv A_{v^{(2)}}\rv - 2m.
\end{align*}
Taking the sum of this inequality over all $v\in V$ and inserting \eqref{eq:v12-against-avearge} gives
\begin{align*}
\sum_{v\in V} \lv A_v\cap I\rv&\geq \sum_{v\in V} \Big( \lv A_v\rv -2m + p^{-(n-3)}\sum_{\substack{w^{(1)}\in V_1, w^{(2)}\in V_2\\w^{(1)}+w^{(2)}=v}} \lv A_{w^{(1)}} \rv + \lv A_{w^{(2)}}\rv\Big)\\
&=\lv A\rv -2mp^{n-1} + p^{-(n-3)}\sum_{w^{(1)}\in V_1,w^{(2)}\in V_2} \lv w^{(1)}+A_{w^{(1)}}(e_1) \rv + \lv w^{(2)}+A_{w^{(2)}}(e_1)\rv\\
&=\lv A\rv -2mp^{n-1} + p\lv A \cap W_1\rv  + p\lv A \cap W_2\rv.
\end{align*}
Combining the above with the bounds $\lv A\rv\geq (m-\alpha)p^{n-1}$ and $\lv A\cap W_i\rv \geq (m-\alpha)p^{n-2}$ for $i\in \{1,2\}$, we find that 
\begin{align*}
\lv A\cap (I(e_1) + V)\rv= \sum_{v\in V} \lv A_v\cap I\rv \geq (m-3\alpha)p^{n-1}.
\end{align*}
We have that $(m-3\alpha)/m>1/2$, whence by a pigeonholing argument, there exists $x\in \F_p$ such $\lv A\cap (\{xe_1,(x+1)e_1\} +V)\rv > p^{n-1}$.
By \Cref{lemma:simple_kneser},
\begin{align*}
e_1+V\subset (A\cap ((x+1)e_1+V))-(A\cap (xe_1+V))
\end{align*}
and
\begin{align*}
-e_1+V\subset (A\cap (xe_1+V))-(A\cap ((x+1)e_1+V)),
\end{align*}
so that $A\cap (\{\pm e_1\} +V)=\emptyset$. Furthermore, we must have that either $\lv A\cap (xe_1 +V)\rv>1/2$ or $\lv A\cap ((x+1)e_1 +V)\rv>1/2$ and hence $A\cap V=\emptyset$, as desired.
\end{proof}

To complete the proof of \Cref{prop:Fp_general}, let $V\subset \F_p^n$ be a subspace of co-dimension 1 such that $A\cap V=A \cap (e_1+V) = A\cap (-e_1+V))=\emptyset$. It remains to show that $V$ satisfies $A\subset I(e_1)+V$.
%\begin{claim}\label{lemma:left-right} If $V\subset \F_p^n$ is a subspace of co-dimension 1 such that $A\cap V=A \cap (e_1+V) = A\cap (-e_1+V))=\emptyset$, then $A\subset I(e_1) + V$.\end{claim}
Let $G$ be the set of good elements inside $V$. Observing that $e_1\notin V$, we see that the family of sets $(v+\langle e_1 \rangle)_{v\in V}$ is a partition of $\F_p^n$, and we can use the bound $\lv A_v\rv \leq m$ to show that $G$ must be large. To be precise, we have
\begin{align*}
(m-\alpha)p^{n-1}\leq \lv A\rv \leq m\lv G\rv + (m-1)(p^{n-1} - \lv G\rv)=mp^{n-1}+\lv G\rv - p^{n-1},
\end{align*}
and thus 
\begin{align}\label{eq:good-lb}
\lv G\rv \geq (1-\alpha)p^{n-1}\geq p^{n-1}/2+p^{n-2},
\end{align} 
which implies that $G+G=V$. Note that by \eqref{eq:offset-addition} the proof is complete if we can show that $t_v=0$ for all good $v$. Indeed, because $G+G=V$, for every $v\in V$ there are elements $v^{(1)},v^{(2)}\in G$ such that $v=v^{(1)}+v^{(2)}$ and thus $A_v\subset t_{v^{(1)}}+t_{v^{(2)}}+I=I$.

In order to achieve this, we first define $L=\{v\in G: t_v\in \{-(m-1),\ldots,0\}\}$, $R=\{v\in G: t_v\in \{0,\ldots,m-1\}$, $C=L\cap R$, $L^+=L\setminus C$, and $R^+=R\setminus C$. Note that $L^+\cup R=L\cup R^+=G$ because $t_v\notin \{-(m-1),\ldots,m-1\}$ would mean that $A_v$ contains 0 and consequently $v\in A\cap V$, contradicting our assumption that $A\cap V=\emptyset$. It is helpful to think of the elements in $C$ as indices to rows that lie in the desired ``central'' position, whereas $L^+$ and $R^+$ indicate rows that lie to the ``left'' or to the ``right'' of this central position. The usefulness of these sets lies in the fact that they are closed under addition in the sense that
\begin{align}\label{eq:left-right-addition}
(L^++ L)\cap G\subset L^+ \qquad\text{and}\qquad(R^++ R)\cap G\subset R^+.
\end{align}
To see this, consider for example $v^{(1)}\in L^+$ and $v^{(2)}\in L$ such that $v^{(1)}+v^{(2)}\in G$. By \eqref{eq:offset-addition}, $t_{v^{(1)}+v^{(2)}}$ is in the set $\{-2(m-1),\ldots,-1\}\subset G\setminus R$ and then in fact in $L^+$ as $G=L^+\cup R$. 

By a similar argument we may also conclude that
\begin{align}\label{eq:left-right-subtraction}
(L^+- R)\cap G\subset L^+ \qquad\text{and}\qquad(R^+- L)\cap G\subset R^+.
\end{align}
We want to show that $C=G$, so let us assume towards a contradiction that at least one of $L^+$ or $R^+$, without loss of generality $L^+$, is non-empty. 

Let us first assume that $R^+$ is empty, so that $L=G$. Denoting by $B=V\setminus G$ the set of bad elements of $V$, we infer from \eqref{eq:left-right-addition} that $L^++L\subset L^+ \cup B$. We know that $L^+$ is a proper subset of $G$ as $0\in C$, whence $L^+ \cup B\neq V$ and $\lv \Sym(L^++L)\rv\leq p^{n-2}$. Consequently, we obtain by Kneser's Theorem that
\begin{align*}
\lv L^+\rv + \lv B\rv \geq \lv L^+ + L \rv \geq \lv L^+ \rv + \lv G\rv - p^{n-2}.
\end{align*}
But this implies $\lv G\rv \leq \lv B\rv + p^{n-2}\leq \alpha p^{n-1}+p^{n-2}=p^{n-1}/2$, contradicting \eqref{eq:good-lb}.

If $R^+$ is not empty, we consider instead the sets $L^+-R$ and $R^+-L$. We know from \eqref{eq:left-right-subtraction} that the former is contained in $L^+\cup B$ and the latter in $R^+\cup B$. Unfortunately, this is not enough to carry out the argument exactly as before, but instead we can make use of the following observation. Whenever we have $v\in (L^+-R)\cap (R^+-L)$, $v$ is not only bad, but \emph{doubly bad}, meaning that $\lv A_v\rv \leq m-2$. Indeed, there must exist $v^{(1)}\in L^+$, $v^{(2)}\in R$, $w^{(1)}\in R^+$, and $w^{(2)}\in L$ such that $v=v^{(1)}-v^{(2)}=w^{(1)}-w^{(2)}$. By \eqref{eq:offset-subtraction}, this implies that
\begin{align}\label{eq:in-intersection}
A_v \subset (t_{v^{(1)}}-t_{v^{(2)}} + I)\cap (t_{w^{(1)}}-t_{w^{(2)}}+I).
\end{align}
We have $s_L:=t_{v^{(1)}}-t_{v^{(2)}}\in \{-2(m-1),\ldots,-1\}$ and $s_R:=t_{w^{(1)}}-t_{w^{(2)}}\in \{1,\ldots,2(m-1)\}$. If $s_L\in \{-2(m-1),\ldots,-m\}$, then $(s_L + I)$ contains at least two of $-1,0$ and $1$. But $A_v$ can contain none of these three points as then $A$ would contain $v-e_1,v$ or $v+e_1$ which is not the case by assumption. In this case, we can thus conclude that $\lv A_v\rv \leq m-2$.

Let us therefore assume that $s_L\in \{-(m-1),\ldots,-1\}$ and likewise that $s_R\in \{1,\ldots,(m-1)\}$. By \eqref{eq:in-intersection}, we have
\begin{align*}
A_v&\subset (\{-(m-1),\ldots,-1\} +I) \cap (\{1,\ldots,(m-1)\}+I)\\
&=\{1,\ldots,2m-2\}\cap \{m+1,\ldots,3m-2\}\\
&=\{m+1,\ldots,2m-2\}.
\end{align*}
The latter set itself has size only $m-2$, whence $\lv A_v\rv\leq m-2$. 

Now that we have established that every $v\in (L^+ - R)\cap (R^+ - L)$ is doubly bad, we may conclude that
\begin{align*}
\lv A \rv =\sum_{v\in V} \lv A_v\rv \leq mp^{n-1} - \lv B\rv - \lv (L^+ - R)\cap (R^+ - L)\rv.
\end{align*}
Combined with $\lv A\rv \geq (m-\alpha)p^{n-1}$, this yields $\lv B\rv + \lv (L^+ - R)\cap (R^+ - L)\rv\leq \alpha p^{n-1}$, and we have by inclusion-exclusion and \eqref{eq:left-right-subtraction} that
\begin{align}\label{eq:LR-inclusion-exclusion}
\lv L^+ - R\rv + \lv R^+ - L\rv
&= \lv (L^+ - R)\cup (R^+ - L)\rv + \lv (L^+ - R)\cap (R^+ - L)\rv\nonumber\\
&= \lv L^+ \cup R^+ \cup B\rv + \lv (L^+ - R)\cap (R^+ - L)\rv\nonumber\\
&= \lv L^+\rv + \lv R^+\rv + \lv B \rv + \lv (L^+ - R)\cap (R^+ - L)\rv\nonumber\\
&< \lv G\rv +\alpha p^{n-1},
\end{align}
where the strictness in the last inequality stems from the fact $0\in C$. For the same reason, we can assert that $\lv \Sym(L^+ - R)\rv , \lv \Sym(R^+ - L)\rv\leq p^{n-2}$, and therefore by Kneser's Theorem,
\begin{align}\label{eq:LR-kneser}
 \lv L^+ - R\rv + \lv R^+ - L\rv \geq 2(\lv G\rv -p^{n-2}).
\end{align}
Combining \eqref{eq:LR-inclusion-exclusion} and \eqref{eq:LR-kneser}, we obtain $\lv G\rv< \alpha p^{n-1}+2p^{n-2}=p^{n-1}/2+p^{n-2}$, which is in contradiction to \eqref{eq:good-lb}. This completes the proof of \Cref{prop:Fp_general}.
\end{proof}
\subsection*{Acknowledgments}
The author is grateful to his supervisor Julia Wolf for her helpful comments on an earlier version of this paper. He would also like to thank Vsevolod Lev for an interesting exchange that led to the discovery of \Cref{example:largest-non-normal}.
\printbibliography
\end{document}